\documentclass[11pt,twoside,final]{article}

\usepackage{amsmath,amsthm,amssymb}
\numberwithin{equation}{section}

\theoremstyle{plain}
\newtheorem{thm}{Theorem}[section]
\newtheorem{prop}[thm]{Proposition}
\newtheorem{lem}[thm]{Lemma}
\newtheorem{cor}[thm]{Corollary}

\theoremstyle{definition}
\newtheorem{prob}{Problem}[section]

\theoremstyle{remark}
\newtheorem{rem}{Remark}[section]

\newcommand{\N}{\mathbb{N}}
\newcommand{\Z}{\mathbb{Z}}
\newcommand{\Q}{\mathbb{Q}}
\newcommand{\R}{\mathbb{R}}
\newcommand{\C}{\mathbb{C}}

\newcommand{\D}{\partial}
\newcommand{\J}{\widetilde{J}}
\newcommand{\w}{\widetilde{w}}
\newcommand{\g}{\widetilde{g}}

\newcommand{\va}{\boldsymbol{a}}
\newcommand{\vb}{\boldsymbol{b}}

\DeclareMathOperator{\ord}{ord}

\newcommand{\deq}{:=}
\newcommand{\HGF}[2]{{}_{#1}F_{#2}}
\newcommand{\HGDop}[2]{{}_{#1}D_{#2}}
\newcommand{\euler}[1][{}]{\mathcal{E}_{#1}}

\newcommand{\Heun}{D_{\mathrm{H}}}

\newcommand{\Ochiai}{D_{\mathrm{O}}}
\newcommand{\Wakayama}{D_{\mathrm{W}}}

\newcommand{\by}[2]{#1, #2}  % \by{family}{personal}
\newcommand{\journal}[1]{{\itshape #1}}
\newcommand{\volume}[1]{{\bfseries #1}}
\begin{document}

\title{Ap\'ery-like numbers arising from
special values of spectral zeta functions for
non-commutative harmonic oscillators}

\author{Kazufumi KIMOTO%
\thanks{Partially supported by Grant-in-Aid for Young 
Scientists  (B) No.16740021.}
{} and Masato WAKAYAMA%
\thanks{Partially supported by Grant-in-Aid for Scientific 
Research (B) No.15340012.}}
\maketitle

\begin{center}
\itshape
Dedicated to Leonhard Euler
on the 299th anniversary of his birthday.
\end{center}

\begin{abstract}
We derive an expression for the value $\zeta_Q(3)$
of the spectral zeta function $\zeta_Q(s)$ studied in \cite{IW1, IW2}
for the non-commutative harmonic oscillator defined in \cite{PW1}
using a Gaussian hypergeometric function.
In this study,
two sequences of rational numbers,
denoted $\J_2(n)$ and $\J_3(n)$,
which can be regarded as analogues
of the \emph{Ap\'ery numbers},
naturally arise and play a key role in obtaining
the expressions for the values $\zeta_Q(2)$ and $\zeta_Q(3)$.
We also show that the numbers
$\J_2(n)$ and $\J_3(n)$ have
congruence relations
like those of the Ap\'ery numbers.

\smallskip

\noindent
{\bfseries Keywords:}
spectral zeta function,
non-commutative harmonic oscillator,
Heun differential equation,
hypergeometric function,
$\zeta(3)$,
Ap\'ery numbers.

\smallskip

\noindent
{\bfseries 2000 Mathematics Subject Classification:}
11M41, 11A07, 33C20.
% 11M41 : Other Dirichlet series and zeta functions
% 11A07 : Congruences; primitive roots; residue systems
% 33C20 : Generalized hypergeometric series, $_pF_q$
\end{abstract}

%====================================================
\section{Introduction}
%====================================================

Let $Q=Q(x,\D_x)$ be the operator 
defined by
\begin{equation}
Q(x,\D_x) \deq
\begin{pmatrix}\alpha & 0 \\ 0 & \beta\end{pmatrix}
\left(-\frac{\D_x^2}2+\frac{x^2}2\right)
+\begin{pmatrix}0 & -1 \\ 1 & 0\end{pmatrix}
\left(x\D_x+\frac12\right),
\end{equation}
where $\D_x \deq \frac{d}{dx}$ and
$\alpha, \beta \in \R$ satisfy $\alpha\beta>1$ \cite{PW1}.
The differential operator $Q$ defines a positive, self-adjoint operator
on $L^2(\R)\otimes\C^2$
with a discrete spectrum
$$
(0<)\lambda_1\le\lambda_2\le\dots\le\lambda_n\le\dots\to\infty.
$$
The operator $Q$
(or a system define by this operator) is called the 
 \emph{non-commutative harmonic oscillator}
(it is sometimes abbreviated NCHO),
because it is considered a generalization
of the ordinary quantum harmonic oscillator
that has an extra non-commutativity
resulting from the non-commutative pair of matrices
$\left(
\left(\begin{smallmatrix}\alpha & 0 \\ 0 & \beta\end{smallmatrix}
\right)
,\,
\left(\begin{smallmatrix}0 & -1 \\ 1 & 0\end{smallmatrix}\right)\right)$.
Actually, when $\alpha=\beta$,
the operator $Q$ is unitarily equivalent to
a pair of harmonic oscillators.

The spectral zeta function $\zeta_Q(s)$ of $Q$ is
defined as the following Dirichlet series \cite{IW1, IW2}:
\begin{equation}
\zeta_Q(s) \deq \sum_{n=1}^\infty \frac1{\lambda_n^s}.
\end{equation}
If $\Re s>1$, then this series converges absolutely,
and hence it defines a holomorphic function
on the half plane $\Re s>1$.
Further, 
$\zeta_Q(s)$ can be meromorphically continued
to the entire plane \cite{IW1}.
We also note that
the function $\zeta_Q(s)$ is regarded as a deformation of
the Riemann zeta function $\zeta(s)$.
In fact, $\zeta_Q(s)=2(2^s-1)\zeta(s)$ when $\alpha=\beta=\sqrt2$, 
because $Q$ is then unitarily equivalent to the operator
$\frac12(-\D_x^2+x^2)
\left(\begin{smallmatrix}1 & 0 \\ 0 & 1\end{smallmatrix}\right)$.

In \cite{IW2},
the first two special values of $\zeta_Q(s)$ are
obtained as follows:
\begin{equation}
\zeta_Q(2)=\frac{(\alpha^{-1}+\beta^{-1})^2}{2(1-\gamma^2)}\left[
(2^2-1)\zeta(2)
+\left(\frac{\alpha^{-1}-\beta^{-1}}{\alpha^{-1}+\beta^{-1}}\right)^{\!2}
\int_0^\infty\!\!\!\!\int_0^\infty\!\! F(t,s;a)dtds
\right],
\end{equation}
\begin{equation}
\zeta_Q(3)=\frac{(\alpha^{-1}+\beta^{-1})^3}{4(1-\gamma^2)^{3/2}}\left[
(2^3-1)\zeta(3)
+3\left(\frac{\alpha^{-1}-\beta^{-1}}{\alpha^{-1}+\beta^{-1}}\right)^{\!2}
\int_0^\infty\!\!\!\!\int_0^\infty\!\!\!\!\int_0^\infty\!\! F(t,s,u;a)dtdsdu
\right],
\end{equation}
where $\gamma\deq(\alpha\beta)^{-1/2}$,
$a\deq(\alpha\beta-1)^{-1/2}$, and
\begin{align}
F(t,s;a) &\deq \frac{e^{-(t+s)/2}}{1-e^{-(t+s)}}
\left\{1+a^2\frac{(1-e^{-2t})(1-e^{-2s})}{(1-e^{-(t+s)})^2}\right\}^{\!-\frac12},\\
F(t,s,u;a) &\deq \frac{e^{-(t+s+u)/2}}{1-e^{-(t+s+u)}}
\left\{1+a^2\frac{(1-e^{-2t})(1-e^{-2(s+u)})}{(1-e^{-(t+s+u)})^2}\right\}^{\!-\frac12}.
\end{align}
We note that involved calculations 
can yield these integral expressions. In addition, the strategy
used there to obtain $\zeta_Q(2)$ and $\zeta_Q(3)$ cannot be 
used for the higher special values
$\zeta_Q(4), \zeta_Q(5), \dots$
(see \cite{IW2} and \S\ref{sec:remarks}).

Using these expressions,
the values $\zeta_Q(2)$ and $\zeta_Q(3)$ are shown to be represented
by a contour integral of a holomorphic solution in the unit disk
of a singly confluent Heun ordinary differential equation
and this equation with an inhomogeneous term, respectively \cite{IW2}.
Based on this result,
in \cite{O}, the following beautiful expression
in terms of
a hypergeometric function
(or a complete elliptic integral)
is obtained for the special value $\zeta_Q(2)$:
\begin{equation}\label{eq:zeta_Q(2)final}
\begin{split}
\zeta_Q(2)&=\frac34\zeta(2)\frac{(\alpha+\beta)^2}{\alpha\beta(\alpha\beta-1)}
\left\{1+\left(\frac{\alpha-\beta}{\alpha+\beta}\,
\HGF21\!\left(\frac14,\frac34;1;\frac1{1-\alpha\beta}\right)
\right)^{\!2}\right\}\\
&=\frac34\zeta(2)\frac{(\alpha+\beta)^2}{\alpha\beta(\alpha\beta-1)}
\left\{1+\left(\frac{\alpha-\beta}{\alpha+\beta}\,
\int_0^{2\pi}\frac{d\theta}{2\pi\sqrt{1+\cos\theta/\sqrt{1-\alpha\beta}}}
\right)^{\!2}\right\}.
\end{split}
\end{equation}
The main purpose of the present paper is to
obtain a similar expression for the special value $\zeta_Q(3)$.
The result is given in the following.
\begin{thm}\label{thm:maintheorem}
We have
\begin{equation}
\begin{split}
\zeta_Q(3)&=
\frac74\zeta(3)
\frac{(\alpha+\beta)^3}{(\alpha\beta(\alpha\beta-1))^{3/2}}\left\{
1+3\left(\frac{\alpha-\beta}{\alpha+\beta}\,
\HGF21\!\left(\frac14,\frac34;1;\frac1{1-\alpha\beta}\right)
\right)^{\!2}\right\}\\
&-\frac{3(\alpha+\beta)(\alpha-\beta)^2}{2\alpha^2\beta^2(\alpha\beta-1)}
\sum_{k=0}^\infty(-1)^k\binom{-\frac12}{k}^{\!3}(\alpha\beta)^{\!-k}
\sum_{j=0}^{k-1}\frac{1}{(2j+1)^3}
\binom{-\frac12}{j}^{\!\!-2}.
\end{split}
\end{equation}
\end{thm}

This paper is organized as follows.
In \S2 we introduce necessary definitions
and present the basic results obtained in \cite{IW2}.
In \S3 we briefly review the result of \cite{O} for $\zeta_Q(2)$.
In \S4 we derive a formula for $\zeta_Q(3)$ in binomial sum form.

Our study of the special values $\zeta_Q(2)$ and $\zeta_Q(3)$
is based on the detailed investigation of
the rational numbers $\J_2(n)$ and $\J_3(n)$ defined in \S2.
More precisely,
the values $\zeta_Q(2)$ and $\zeta_Q(3)$ are obtained
as series of the numbers $\J_2(n)$ and $\J_3(n)$,
and the calculations of $\zeta_Q(2)$ and $\zeta_Q(3)$ are reduced
to those of $\J_2(n)$ and $\J_3(n)$.
The numbers $\J_2(n)$ and $\J_3(n)$ can be regarded
as analogues of the \emph{Ap\'ery numbers}
(see Remark \ref{rem:aperylike}),
which have been used to prove
the irrationality of $\zeta(2)=\frac{\pi^2}6$ and $\zeta(3)$
(see, e.g., \cite{Po1979}; see also \cite{AAR}). 
It is quite mysterious that
although these proofs of irrationality
for $\zeta(2)$ and $\zeta(3)$ are parallel,
the same method cannot be used to prove
the irrationality of the higher special values
$\zeta(4), \zeta(5), \dots,$ etc.
In our spectral zeta function case,
the situation is similar in the sense that
the special values $\zeta_Q(2)$ and $\zeta_Q(3)$ can be calculated
and described in quite similar ways
using the numbers $\J_2(n)$ and $\J_n(3)$ (see \S2),
but the same method (and/or strategy) does not seem to be effective for
calculations of the higher special values
$\zeta_Q(4), \zeta_Q(5), \dots,$ etc.
Further, the generating functions for $\J_2(n)$ and $\J_3(n)$
satisfy certain
(singly confluent Heun) differential equations in the same operator,
with that for $\zeta_Q(2)$ being homogeneous and
that for $\zeta_Q(3)$ being inhomogeneous
(see \eqref{de_for_w2} and \eqref{de_for_w3}) \cite{IW2}.
It is important to investigate
the cause of those special situations
for the pairs $\zeta(2)$ and $\zeta(3)$
and $\zeta_Q(2)$ and $\zeta_Q(3)$.
We expect that
there might be some relationship between $\zeta_Q(2)$ and $\zeta_Q(3)$
similar to Euler's conjecture \cite{Euler1785},
\begin{equation}\label{eq:Euler}
\zeta(3)=\alpha(\log2)^3+\beta\zeta(2)\log2
\qquad(\exists\alpha,\exists\beta\in\Q).
\end{equation}

In \S\ref{sec:remarks}, we give some remarks and
introduce some problems
concerning the integrals which might be 
related to the higher special values $\zeta_Q(k)\,(k\geq4)$. 
As an important result of the recent studies
of the irrationality of special values
of the Riemann zeta function
that differs from the conventional Ap\'ery-Beukers line of study,
we refer the reader to Rivoal \cite{R2002}.

In addition to our main purpose,
there are several points of interest in the present study.
First, we note that
there are a number of interesting studies
of the Ap\'ery numbers
concerning congruence properties
and their analogues
(see, e.g.,
\cite{CCC1980,Co1980,Ge1982,Be1985,Be1987,AO2000,OP,Ko1992,Mi1983}
and references therein).
These studies suggest that
the congruence originates from
automorphic forms associated with some algebraic surface
(see, e.g., \cite{BP1984,SB1985}).
In \S\ref{sec:congruence} we study several basic
(Ap\'ery number-like)
congruence properties of the numbers $\J_2(n)$ and $\J_3(n)$.
It would be interesting to investigate the interpretation of
$\zeta_Q(2)$ and $\zeta_Q(3)$ as period integrals
on the algebraic surfaces
from which the congruence relations among
the numbers $\J_2(n)$ and $\J_3(n)$ arise.

In the appendix, \S\ref{sec:appendix},
we generalize the method used in \S4
and derive a general formula for
the holomorphic solutions of inhomogeneous hypergeometric
differential equations.
As an application,
we present another method of deriving an explicit expression
for the numbers $\J_3(n)$.

Throughout the paper, we denote by
$\N,\Z,\Q,\R$ and $\C$ the set of natural numbers,
rational integers, rational numbers, real numbers
and complex numbers, respectively.
Also, we denote by $\HGF{p}{q}(a_1,\dots,a_p;b_1,\dots,b_q;z)$
the (generalized) hypergeometric function (see \S\ref{sec:appendix}).

%====================================================
\section{Definitions and basic results}\label{sec:preliminaries}
%====================================================

We first recall the necessary definitions
and conventions from \cite{IW2}
and present the basic properties.
(For detailed discussions, see \cite{IW2}.)

If $\alpha\beta>2$ or $a^2=(\alpha\beta-1)^{-1}<1$,
the functions $F(t,s;a)$ and $F(t,s,u;a)$ have Taylor expansions
with respect to $a^2$ around the origin,
because
$(1-e^{-2t})(1-e^{-2(s+u)})\le(1-e^{-(t+s+u)})^2$
in the case $t,s,u>0$.
Therefore, with the definitions
\begin{align}
J_2(n) &\deq
\int_0^\infty\!\!\!\!\int_0^\infty\!\!
\frac{e^{-(t+s)/2}}{1-e^{-(t+s)}}
\left(\frac{(1-e^{-2t})(1-e^{-2s})}{(1-e^{-(t+s)})^2}\right)^{\!n}
dtds,\label{eq:integral_exp_of_J2}\\
J_3(n) &\deq
\int_0^\infty\!\!\!\!\int_0^\infty\!\!\!\!\int_0^\infty\!\!
\frac{e^{-(t+s+u)/2}}{1-e^{-(t+s+u)}}
\left(\frac{(1-e^{-2t})(1-e^{-2(s+u)})}{(1-e^{-(t+s+u)})^2}\right)^{\!n}
dtdsdu\label{eq:integral_exp_of_J3}
\end{align}
and
\begin{align}
g_2(x)&\deq\sum_{n=0}^\infty\binom{-\frac12}{n}J_2(n)x^n,\\
g_3(x)&\deq\sum_{n=0}^\infty\binom{-\frac12}{n}J_3(n)x^n,
\end{align}
the special values $\zeta_Q(2)$ and $\zeta_Q(3)$ can be expressed as
\begin{align}
\zeta_Q(2) &= 
\frac{(\alpha^{-1}+\beta^{-1})^2}{2(1-\gamma^2)}\left[
(2^2-1)\zeta(2)+\left(\frac{\alpha-\beta}{\alpha+\beta}\right)^2
g_2(a^2)\right],\label{eq:zetaQ(2)}\\
\zeta_Q(3) &= \frac{(\alpha^{-1}+\beta^{-1})^3}{4(1-\gamma^2)^{3/2}}\left[
(2^3-1)\zeta(3)
+3\left(\frac{\alpha-\beta}{\alpha+\beta}\right)^2
g_3(a^2)\label{eq:zetaQ(3)}
\right]
\end{align}
if $\alpha\beta>2$.
From this point, we assume that $\alpha\beta>2$.

The numbers $J_k(n)$ satisfy the following recurrence relations.

\begin{prop}[{\cite[Propositions 4.11 and 6.4]{IW2}}]
We have
\begin{align}
4n^2 J_2(n)-(8n^2-8n+3)J_2(n-1)+4(n-1)^2J_2(n-2) &= 0,\label{rec_for_J2}\\
4n^2 J_3(n)-(8n^2-8n+3)J_3(n-1)+4(n-1)^2J_3(n-2) &= \frac{2^n(n-1)!}{(2n-1)!!},\label{rec_for_J3}
\end{align}
with the initial values
\begin{align*}
&J_2(0)=3\zeta(2),\quad J_2(1)=\frac94\zeta(2),\\
&J_3(0)=7\zeta(3),\quad J_3(1)=\frac{21}4\zeta(3)+\frac12.
\end{align*}
\qed
\end{prop}

\begin{rem}\label{rem:aperylike}
Roger Ap\'ery introduced the numbers
$A_n=\sum_{k=0}^n\binom{n}{k}^2\binom{n+k}{k}^2$
in 1978
called \emph{Ap\'ery numbers} nowadays,
which satisfy the recurrence relation
\begin{equation}\label{eq:apery_recurrence}
A_0=1,A_1=5,\quad
(n+1)^3A_{n+1}-(34n^3+51n^2+27n+5)A_n+n^3A_n=0
\quad(n\ge1).
\end{equation}
Ap\'ery used these numbers 
to prove the irrationality of $\zeta(3)=\sum_{n=1}^\infty n^{-3}$
(see, e.g., \cite{Po1979}). Note that
the recurrence relations
\eqref{rec_for_J2} and \eqref{rec_for_J3} for
$\J_2(n)$ and $\J_3(n)$
are quite similar to \eqref{eq:apery_recurrence}.
Moreover,
$\J_2(n)$ and $\J_3(n)$ satisfy congruence relations
analogous to those satisfied by the Ap\'ery numbers $A_n$
(see \S\ref{sec:congruence}).
In this sense,
we may regard $\J_n(n)$ and $\J_3(n)$ as analogues of
the Ap\'ery numbers $A_n$.
\end{rem}

Denote by $w_k(z)$ the generating functions
for the numbers $J_k(n)$, defined by
\begin{equation}
w_k(z) \deq \sum_{n=0}^\infty J_k(n)z^n.
\end{equation}
From the recurrence relations \eqref{rec_for_J2} and \eqref{rec_for_J3},
it follows that each of the functions $w_2(z)$ and $w_3(z)$
satisfies a confluent Heun differential equation,
as described by the following

\begin{lem}[{\cite[Theorems 4.13 and 6.5]{IW2}}]
The functions $w_2(z)$ and $w_3(z)$ satisfy
\begin{align}
\Heun w_2(z) &= 0,\label{de_for_w2}\\
\Heun w_3(z) &= \frac12 \HGF21\!\left(1,1;\frac32;z\right),\label{de_for_w3}
\end{align}
where $D_H$ is the singly confluent
Heun differential operator defined by
\begin{equation}
\Heun \deq z(1-z)^2\D_z^2+(1-3z)(1-z)\D_z+z-\frac34.
\end{equation}
\qed
\end{lem}

For convenience, we introduce the `normalized' sequences
\begin{equation}
\J_2(n)\deq \frac{J_2(n)}{J_2(0)},\qquad
\J_3(n)\deq J_3(n)-J_3(0)\J_2(n).
\end{equation}
These numbers also satisfy the recurrence relations
\eqref{rec_for_J2} and \eqref{rec_for_J3},
with the initial values
\begin{align*}
&\J_2(0)=1,\quad \J_2(1)=\frac34,\\
&\J_3(0)=0,\quad \J_3(1)=\frac12.
\end{align*}
In particular, it is easy to see that $\J_2(n), \J_3(n)\in\Q$.
Numerically, we have

\begin{table}[htb]\label{Table1}
\begin{center}
\begin{tabular}{c|cccccccccc}
$n$ & $0$ & $1$ & $2$ & $3$ & $4$ & $5$ & $6$ & $7$ & $8$ & $9$ \\
\hline
\rule[-0.5ex]{0ex}{4.5ex}
$\J_2(n)$ &
$1$ & $\dfrac{3}{4}$ & $\dfrac{41}{64}$ & $\dfrac{147}{256}$ & $\dfrac{8649}{16384}$ &
$\dfrac{32307}{65536}$ & $\dfrac{487889}{1048576}$ & $\dfrac{1856307}{4194304}$ &
$\dfrac{454689481}{1073741824}$ & $\dfrac{1748274987}{4294967296}$
\end{tabular}
\caption{First ten values of $\J_2(n)=J_2(n)/J_2(0)$.}
\end{center}
\end{table}

\begin{table}[htb]\label{Table2}
\begin{center}
\begin{tabular}{c|ccccccc}
$n$ & $0$ & $1$ & $2$ & $3$ & $4$ & $5$ & $6$ \\
\hline
\rule[-0.5ex]{0ex}{4.5ex}
$\J_3(n)$ &
$0$ & $\dfrac{1}{2}$ & $\dfrac{65}{96}$ & $\dfrac{13247}{17280}$ & $\dfrac{704707}{860160}$ & $\dfrac{660278641}{774144000}$ & $\dfrac{357852111131}{408748032000}$
\end{tabular}

\medskip

\begin{tabular}{c|ccc}
$n$ & $7$ & $8$ & $9$ \\
\hline
\rule[-0.5ex]{0ex}{4.5ex}
$\J_3(n)$ &
$\dfrac{309349386395887}{347163328512000}$ & $\dfrac{240498440880062263}{266621436297216000}$ & $\dfrac{148443546307725010253}{163172319013896192000}$
\end{tabular}
\caption{First ten values of $\J_3(n)=J_3(n)-J_3(0)\J_2(n)$.}
\end{center}
\end{table}
We note that the convergence of $\J_3(n)$ to $0$ is slower than
that of $\J_2(n)$
(see Table \ref{Table2} below).
For instance, we have $\J_2(10^4)\cong 0.025$ and $\J_3(10^4)\cong 0.2457$.
%We note that the convergence of $\J_3(n)$ to $0$ is quite moderate
%(see Table \ref{Table2} below).
%For instance, we have $\J_3(10^4)\cong 0.2457$.
% $\J_3(n) > \J_2(n)\log n$ when $n \gg 1$ ?

\medskip

We next introduce the generating functions
of $\J_2(n)$ and $\J_3(n)$,
$\w_2(x)$, $\w_3(x)$, $\g_2(x)$ and $\g_3(x)$,
in analogy to $w_k(z)$ and $g_k(z)$:
\begin{align}
\w_k(x)&\deq \sum_{n=0}^\infty \J_k(n)x^n,\label{eq:def_for_wt_k}\\
\g_k(x)&\deq \sum_{n=0}^\infty \binom{-\frac12}{n}\J_k(n)x^n.\label{eq:def_for_gt_k}
\end{align}
It is clear that
the function $\w_k(z)$ satisfies
the same differential equation as $w_k(z)$.

%====================================================
\section{Explicit formulas for $\J_2(n)$ and $\zeta_Q(2)$}\label{sec:J2}
%====================================================

We briefly review the derivation of the formula for $\zeta_Q(2)$
according to Ochiai \cite{O},
and give an expression for $\J_2(n)$
in terms of binomial coefficients as its corollary.

We first see that the function $w_2(z)$ is expressed
in terms of a hypergeometric function as follows.
\begin{prop}[{\cite[Proposition 3]{O}}]
The generating function $w_2(z)$ is given by
\begin{equation}
w_2(z)=
\frac{J_2(0)}{1-z}\,\HGF21\!\left(\frac12,\frac12;1;\frac{z}{z-1}\right).
\end{equation}
\end{prop}

\begin{proof}
In order to make the paper self-contained,
we present here the proof given by Ochiai \cite{O}.
It is straightforward to verify that the equation 
\eqref{de_for_w2} is equivalent to
\begin{equation}\label{eq:OchiaiReduction}
4(1-z)^2\D_z z\D_z(1-z)w_2(z)+(1-z)w_2(z)=0.
\end{equation}
Then, changing the variable by $t=\frac{z}{z-1}$
and putting $v_2(t)=(1-z)w_2(z)$,
we see that \eqref{eq:OchiaiReduction} is equivalent to
the hypergeometric differential equation
\begin{equation}
\left(t(1-t)\D_t^2+(1-2t)\D_t-\frac14\right)v_2(t)=0,
\end{equation}
whose local holomorphic solution is a scalar multiple of
$\HGF21(\frac12,\frac12;1;t)$.
This completes the proof.
\end{proof}

As a corollary,
we have the following explicit expression for $\J_2(n)$.
\begin{cor}
The number $\J_2(n)$ is given by
\begin{equation}\label{eq:explicit_expression_for_J2}
\J_2(n)=\sum_{k=0}^n (-1)^k \binom{-\frac12}{k}^{\!\!2}\binom{n}{k}.
\end{equation}
\qed
\end{cor}

The expression
\eqref{eq:zeta_Q(2)final} for $\zeta_Q(2)$
in terms of the hypergeometric function
is obtained by combining \eqref{eq:zetaQ(2)} and
the formula \eqref{g2formula_final} for $\g_2(x)$,
proved below.
The proof is essentially the same as that given by Ochiai \cite{O},
however,
in contrast to the proof of Ochiai,
our proof does not employ integrations.

\begin{prop}
We have
\begin{equation}\label{g2formula_final}
\begin{split}
\g_2(x)
=\frac1{\sqrt{1+x}}\HGF21\!\left(\frac14,\frac14;1;\frac{x}{1+x}\right)^{\!2}
=\HGF21\!\left(\frac14,\frac34;1;-x\right)^{\!2}.
\end{split}
\end{equation}
\end{prop}

\begin{proof}
Recall the following two elementary formulas for binomial sums:
\begin{align}
\sum_{n=0}^\infty(-1)^n\binom{-\frac12}{n}^{\!\!3}x^n
=\HGF32\!\left(\frac12,\frac12,\frac12;1,1;x\right),\\
\sum_{n=k}^\infty \binom{-\frac12}{n}\binom{n}{k}x^n
= \frac1{\sqrt{1+x}}\binom{-\frac12}{k}\left(\frac{x}{1+x}\right)^{\!\!k}.
\end{align}
Using these, we can calculate
\begin{equation}\label{eq:calc_of_g2}
\begin{split}
\g_2(x)&=\sum_{n=0}^\infty\binom{-\frac12}{n}x^n
\sum_{k=0}^n (-1)^k \binom{-\frac12}{k}^{\!\!2}\binom{n}{k}\\
&=\sum_{k=0}^\infty (-1)^k \binom{-\frac12}{k}^{\!\!2}
\sum_{n=k}^\infty \binom{-\frac12}{n}\binom{n}{k}x^n\\
&=\frac1{\sqrt{1+x}}
\sum_{k=0}^\infty (-1)^k \binom{-\frac12}{k}^{\!\!3}
\left(\frac{x}{1+x}\right)^{\!\!k}\\
&=\frac1{\sqrt{1+x}}\HGF32\!\left(\frac12,\frac12,\frac12;1,1;\frac{x}{1+x}\right).
\end{split}
\end{equation}
We then apply the formulas (see, e.g., \cite{AAR})
\begin{align}
\HGF21\!\left(\alpha,\beta;\alpha+\beta+\frac12;z\right)^{\!2}&=
\HGF32\!\left(2\alpha,2\beta,\alpha+\beta;2\alpha+2\beta,\alpha+\beta+\frac12;z\right),
\tag{Clausen's identity}\\
\HGF21(\alpha,\beta;\gamma;z)&=
(1-z)^{-\alpha}\HGF21\!\left(\alpha,\gamma-\beta;\gamma;\frac{z}{z-1}\right)
\tag{Pfaff's formula}
\end{align}
to \eqref{eq:calc_of_g2}.
This gives the desired result.
\end{proof}

\begin{rem}
We immediately deduce the following expression for $\J_2(n)$
from \eqref{g2formula_final}:
\begin{equation}
\J_2(n)=2^{-4n}\binom{2n}{n}^{-1}
\sum_{k=0}^n \binom{2k}{k}\binom{4k}{2k}
\binom{2n-2k}{n-k}\binom{4n-4k}{2n-2k}.
\end{equation}
This clearly implies the positivity of $\J_2(n)$.
\end{rem}

%====================================================
\section{Binomial sum formulas for $\J_3(n)$ and $\zeta_Q(3)$}\label{sec:J3}
%====================================================

In this section,
we present a binomial sum formula for $\J_3(n)$
like \eqref{eq:explicit_expression_for_J2} for $\J_2(n)$,
and give an expression for $\zeta_Q(3)$
(Theorem \ref{thm:maintheorem_revisited})
like \eqref{eq:zeta_Q(2)final} for $\zeta_Q(2)$.
The strategy for obtaining these is the same as that used in \S3.
Specifically, we solve the differential equation \eqref{de_for_w3}
and utilize the solution to obtain a formula for $\J_3(n)$.

We construct a holomorphic solution $w(z)$ of the differential equation
\begin{equation}\label{H1}
\Heun w(z)=\frac12\HGF21\!\left(1,1;\frac32;z\right)
\end{equation}
at the origin. For this purpose,
first, we change the variable by $t=\dfrac{z}{z-1}$
and introduce
\begin{align*}
\Ochiai&\deq t(1-t)\D_t^2+(1-2t)\D_t-\frac14,\\
v(t)&\deq(1-z)w(z).
\end{align*}
Then, we see that the equation \eqref{H1} is equivalent to
\begin{equation}\label{G1}
\Ochiai v(t)=-\frac12\frac1{1-t}\HGF21\!\left(1,1;\frac32;\frac{t}{t-1}\right).
\end{equation}
The right-hand side of \eqref{G1} takes the simple form given in the following.
\begin{lem}
We have
\begin{equation}
\frac1{1-t}\HGF21\!\left(1,1;\frac32;\frac{t}{t-1}\right)
=\sum_{n=0}^\infty\frac{t^n}{2n+1}.
\end{equation}
\end{lem}

\begin{proof}
Applying Pfaff's formula, we find
\begin{equation*}
\frac1{1-t}\HGF21\!\left(1,1;\frac32;\frac{t}{t-1}\right)
=\HGF21\!\left(1,\frac12;\frac32;t\right)
=\sum_{n=0}^\infty\frac{t^n}{2n+1}.
\end{equation*}
\end{proof}

\begin{lem}
The polynomial function $p_n(t)$ defined as
\begin{equation}
p_n(t)=-\frac{4}{(2n+1)^2}\binom{-\frac12}{n}^{\!\!\!-2}
\sum_{k=0}^n \binom{-\frac12}{k}^{\!\!2} t^k
\end{equation}
satisfies the differential equation $\Ochiai p_n(t)=t^n$.
\end{lem}

\begin{proof}
The assertion is verified by
straightforward calculation.
(See also Proposition \ref{prop:polynomial}.)
\end{proof}

\medskip

Thus, the function
\begin{equation}
v(t)\deq -\frac12\sum_{n=0}^\infty\frac{p_n(t)}{2n+1}
=2\sum_{k=0}^\infty\left(
\sum_{n=k}^\infty \frac{1}{(2n+1)^3}\binom{-\frac12}{n}^{\!\!\!-2}
\right)\binom{-\frac12}{k}^{\!\!2}t^k
\end{equation}
gives a local holomorphic solution to \eqref{G1}.
Therefore,
\begin{equation}
w(z)=\frac1{1-z}v\left(\frac{z}{z-1}\right)
=2\sum_{n=0}^\infty
\left(
\sum_{k=0}^n(-1)^k\binom{n}{k}\binom{-\frac12}{k}^{\!\!2}
\sum_{j=k}^\infty \frac{1}{(2j+1)^3}\binom{-\frac12}{j}^{\!\!\!-2}
\right)z^n
\end{equation}
is a local holomorphic solution to \eqref{H1}.

Note that any holomorphic solution $w$ of \eqref{H1} at the origin
can be written in the form $w=c\w_2+\w_3$ for some constant $c$.
Then, it is seen that
the constant $c$ is simply the constant term of $w$;
that is, we have $\w_3=w-w(0)\w_2$.
The constant term of $w$ is given by
$$
w(0)=v(0)=2\sum_{j=0}^\infty \frac{1}{(2j+1)^3}\binom{-\frac12}{j}^{\!\!\!-2}.
$$
Thus we have the

\begin{thm}\label{thm:J3-formula}
The holomorphic solution $\w_3$ of the differential equation \eqref{H1}
at the origin with the initial condition
$\w_3(0)=0$ is given by
\begin{equation}\label{w3}
\begin{split}
\w_3(z)
&=-2\sum_{n=0}^\infty
\left(
\sum_{k=0}^n(-1)^k\binom{-\frac12}{k}^{\!\!2}\binom{n}{k}
\sum_{j=0}^{k-1} \frac{1}{(2j+1)^3}\binom{-\frac12}{j}^{\!\!\!-2}
\right)z^n.
\end{split}
\end{equation}
In particular, 
\begin{equation}\label{eq:formula_for_J3}
\J_3(n)=-2\sum_{k=0}^n(-1)^k\binom{-\frac12}{k}^{\!\!2}\binom{n}{k}
\sum_{j=0}^{k-1} \frac{1}{(2j+1)^3}\binom{-\frac12}{j}^{\!\!\!-2}.
\end{equation}
\qed
\end{thm}

With the above, the function $\g_3(x)$ can be calculated
in the same manner as $\g_2(z)$
in the previous section.

\begin{cor}
The function $\g_3(x)$ is given by
\begin{equation}\label{eq:g3}
\g_3(x)=\frac{-2}{\sqrt{1+x}}
\sum_{k=1}^\infty(-1)^k\binom{-\frac12}{k}^{\!\!3}\left(\frac{x}{1+x}\right)^{\!\!k}
\sum_{j=0}^{k-1}\frac1{(2j+1)^3}\binom{-\frac12}{j}^{\!\!\!-2}.
\end{equation}
\qed
\end{cor}

Combining \eqref{eq:g3} and \eqref{eq:zetaQ(3)},
we obtain the
\begin{thm}[Theorem \ref{thm:maintheorem}]\label{thm:maintheorem_revisited}
The value $\zeta_Q(3)$ is given by
\begin{equation}
\begin{split}
\zeta_Q(3)&=
\frac74\zeta(3)
\frac{(\alpha+\beta)^3}{(\alpha\beta(\alpha\beta-1))^{3/2}}\left\{
1+3\left(\frac{\alpha-\beta}{\alpha+\beta}\,
\HGF21\!\left(\frac14,\frac34;1;\frac1{1-\alpha\beta}\right)
\right)^{\!2}\right\}\\
&-\frac{3(\alpha+\beta)(\alpha-\beta)^2}{2\alpha^2\beta^2(\alpha\beta-1)}
\sum_{k=0}^\infty\binom{-\frac12}{k}^{\!3}\left(\frac{-1}{\alpha\beta}\right)^{\!k}
\sum_{j=0}^{k-1}\frac{1}{(2j+1)^3}
\binom{-\frac12}{j}^{\!\!-2}.
\end{split}
\end{equation}
\qed
\end{thm}

%====================================================
\section{Remarks related to higher special values}\label{sec:remarks}
%====================================================

Noting the integral expressions
\eqref{eq:integral_exp_of_J2} and
\eqref{eq:integral_exp_of_J3}
for $J_2(n)$ and $J_3(n)$,
it is quite natural to introduce the numbers $J_k(n)$ through the integral
\begin{align}
J_k(n) &\deq
\int_0^\infty\!\!\!\!\int_0^\infty\!\!\!\!\dotsi\!\int_0^\infty\!\!
\frac{e^{-(t_1+\dotsb+t_k)/2}}{1-e^{-(t_1+\dotsb+t_k)}}
\left(\frac{(1-e^{-2t_1})(1-e^{-2(t_2+\dotsb+t_k)})}
{(1-e^{-(t_1+\dotsb+t_k)})^2}\right)^n
dt_1dt_2\dotsb dt_k\\
&=2^k\int_0^1\!\!\int_0^1\!\!\dotsi\!\int_0^1\!\!
\frac{1}{1-x_1^2\dotsb x_k^2}
\left(\frac{(1-x_1^4)(1-x_2^4\dotsb x_k^4)}{(1-x_1^2\dotsb x_k^2)^2}\right)^n
dx_1dx_2\dotsb dx_k.
\end{align}
Then, for their generating functions, we have
\begin{equation}
w_k(z) \deq \sum_{n=0}^\infty J_k(n)z^n
=2^k\int_0^1\!\!\int_0^1\!\!\dotsi\!\int_0^1\!\!
\frac{1-x_1^2\dotsb x_k^2}{(1-x_1^2\dotsb x_k^2)^2-(1-x_1^4)(1-x_2^4\dotsb x_k^4)z}
dx_1dx_2\dotsb dx_k.
\end{equation}
Note that the numbers $J_k(n)$ are all positive. 
Also we have $w_k(0)=(2^k-1)\zeta(k)$.

\begin{rem}
Because in general,
for a sufficiently well-behaved $f$, we have the relation
\begin{equation}
\int_0^\infty\!\!\!\!\int_0^\infty\!\!\!\!\dotsi\!\int_0^\infty\!\!
f(x_1+\dots+x_n)dx_1\dotsb dx_n=
\frac1{\Gamma(n)}\int_0^\infty y^{n-1}f(y)dy,
\end{equation}
we see that
\begin{equation}
J_k(n) =
\frac1{\Gamma(k-1)}
\int_0^\infty\!\!\int_0^\infty
s^{k-2}\frac{e^{-(s+t)/2}}{1-e^{-(s+t)}}
\left(\frac{(1-e^{-2s})(1-e^{-2t})}
{(1-e^{-(s+t)})^2}\right)^n
dsdt.
\end{equation}
\end{rem}

Here, we point out that
it is still unclear if the special values 
$\zeta_Q(k)\,(k\geq4)$
of the spectral zeta function for the NCHO 
can also be described by the numbers $J_k(n)$.
More precisely, we have the following question: 
%\begin{quest}
Is it true that
\begin{equation}
\zeta_Q(k) = Z_0(k)+\sum_{n=0}^\infty Z'_n(k)
\qquad(k\ge2)
\end{equation}
with
\begin{align}
Z_0(k)&=\frac{(\alpha^{-1}+\beta^{-1})^k}{2^{k-1}(1-\gamma^2)^{k/2}}
(2^k-1)\zeta(k),\\
Z'_n(k)&=\frac{(\alpha^{-1}+\beta^{-1})^k}{2^{k-1}(1-\gamma^2)^{k/2}}
R_{n,k}\!\left(\frac{\alpha-\beta}{\alpha+\beta}\right)\!
(-1)^n\binom{2n}{n}\left(\frac{a}2\right)^{2n}J_k(n)
\end{align}
for some polynomial $R_{n,k}$ ?
%\qed
%\end{quest}
\medskip

Beside this question, we have the 
\begin{prob}
Find a differential equation satisfied by $w_k(z)\,(k\geq4)$. \qed
\end{prob}

Now, note that
the numbers $J_k(1)$ can be written
\begin{equation}
J_k(1) = \frac34\left(
\zeta(k,\frac12)+\sum_{m=1}^{[k/2]-1}2^{-2m}\zeta(k-2m,\frac12)
\right)+\frac{1-(-1)^k}{2^{k-1}},
\end{equation}
where $\zeta(s,x)=\sum_{n=0}^\infty (n+x)^{-s}$ is the Hurwitz zeta function.
This implies that $J_{k-2}(1)-4J_k(1)+3J_k(0)=0$.
We are thus led to the following.

\begin{prob}
Carry out a general study of
such ``vertical'' relations among the numbers $J_k(n)$. \qed
\end{prob}

%====================================================
\section{Congruence properties of $\J_2(n)$ and $\J_3(n)$}\label{sec:congruence}
%====================================================

The numbers $\J_2(n)$ and $\J_3(n)$, as we have seen, can be regarded
as analogues of the Ap\'ery numbers,
$A_n=\sum_{k=0}^n\binom{n}{k}^2\binom{n+k}{k}^2$ (see Remark 
\ref{rem:aperylike}).
In this section we study congruence properties
of $\J_2(n)$ and $\J_3(n)$ similar to those of
the Ap\'ery numbers.

We employ the following convention for treating
congruence properties of rational numbers.
Let $x, y \in \Q$ and $r$ be a positive integer.
If $x-y = \frac{rn}{d}$ for some $n, d \in \Z$
such that $d$ is relatively prime to $r$,
then we regard $x$ and $y$ are congruent modulo $r$
and write $x \equiv y \pmod{r}$.

We now demonstrate the following basic congruence
(or divisibility) property of $\J_2(n)$.
\begin{prop}
Let $p$ be an odd prime number such that $p\equiv-1\pmod4$.
We have $\J_2(n)\equiv0\pmod p$ if
$n_j(p)\equiv\frac{p-1}2$ for some $j$,
where $n=n_0(p)+n_1(p)p+\dotsb+n_d(p)p^d$ is the base $p$ expansion of $n$.
\end{prop}

\begin{proof}
First, note that
\begin{equation}\label{eq:base-p_reduction}
\J_2(n) \equiv \J_2(n_0(p))\J_2(n_1(p))\dotsb\J_2(n_d(p)) \pmod{p}.
\end{equation}
This relation can be obtained by repeatedly using
the elementary congruence formulas
\begin{equation}
\binom{ap+b}{cp+d} \equiv \binom{a}{c}\binom{b}{d} \pmod{p},\qquad
2^{-p} \equiv 2^{-1} \pmod{p}
\end{equation}
in \eqref{eq:explicit_expression_for_J2},
the explicit expression for $\J_2(n)$.

Next, we show that
$\J_2\left(\frac{p-1}2\right)\equiv0$
if $p\equiv-1\pmod4$.
Because
$\binom{\frac{p-1}2}{k}
\equiv(-1)^k 2^{-2k} \binom{2k}{k}
\pmod{p}$,
it follows from \eqref{eq:explicit_expression_for_J2} that
\begin{equation*}
\begin{split}
\J_2\left(\frac{p-1}2\right)
\equiv 
\sum_{k=0}^{\frac{p-1}2}
(-1)^k\binom{\frac{p-1}2}{k}^{\!\!3}
\pmod{p}.
\end{split}
\end{equation*}
Next, using the formulas
\begin{equation*}
\sum_{k=0}^n (-1)^k \binom{n}{k}^{\!3} = \begin{cases}
\,0\, & \text{$n$\,:\,odd},\\
(-1)^{n/2}\frac{(3n/2)!}{(n/2)!^3} & \text{$n$\,:\,even},
\end{cases}
\end{equation*}
\begin{equation*}
\left(\frac{3(p-1)}4\right)!\equiv
(-1)^{(p-1)/4}{\left(\frac{p-1}4\right)!}^{\!\!-1} \pmod{p},
\end{equation*}
we have
\begin{equation}\label{eq:J_2(p-1/2)}
\J_2\left(\frac{p-1}2\right)\equiv
\begin{cases}
\,0\, & p \equiv -1 \pmod4,\\
\displaystyle
-{\left(\tfrac{p-1}4\right)!}^{\!\!-4} & p \equiv 1 \pmod4,
\end{cases}
\end{equation}
as desired.
Combining \eqref{eq:base-p_reduction} and \eqref{eq:J_2(p-1/2)},
we have the conclusion.
\end{proof}

\begin{rem}
The condition that
$n_j(p)\equiv\frac{p-1}2$ holds for some $j$
is not a necessary condition
to realize $\J_2(n)\equiv0\pmod{p}$.
In fact, for a given prime number $p$ such that $p\equiv3\pmod4$,
there can exist $n\in\N$ distinct from $(p-1)/2$ satisfying
$\J_2(n) \equiv 0 \pmod{p}$, with $0 \le n < p$.
For instance, we have
$\J_2(7) \equiv \J_2(15) \equiv 0 \pmod{23}$.
We also note that
there exists a prime number $p$ which is congruent to $1$ modulo $4$
such that $\J_2(n) \equiv 0 \pmod{p}$ for some $n$ with $0 \le n < p$.
\end{rem}

We now show higher-order congruence properties of $\J_2(n)$ and $\J_3(n)$.
Each of these is an analogue of the congruence
\begin{equation}
A_{mp^n-1} \equiv A_{mp^{n-1}-1} \pmod{p^{n}}
\end{equation}
for Ap\'ery numbers.

\begin{thm}\label{prop:J2_higher_congruence}
For any odd prime number $p$, we have
\begin{equation}\label{eq:higher_congruence_2}
\J_2(mp^r) \equiv \J_2(mp^{r-1}) \pmod{p^r}
\end{equation}
\begin{equation}\label{eq:higher_congruence_3}
\J_3(p^r)p^{3r} \equiv \J_3(p^{r-1})p^{3(r-1)} \pmod{p^r}
\end{equation}
for any $m,r\in\N$.
\end{thm}

\begin{proof}
First, note the elementary congruence relation
\begin{equation}\label{eq:basic_binom_congruence}
\binom{mp^r}{kp^l}\equiv\begin{cases}
\binom{mp^{r-1}}{kp^{l-1}} & l\ge1\\
0 & l=0
\end{cases}
\pmod{p^r}
\end{equation}
for $k \not\equiv 0 \pmod{p}$,
from which it follows immediately that
\begin{equation}\label{eq:needed-1}
\binom{mp^r}{kp^l}\equiv0\pmod{p^{r-l}}.
\end{equation}
We also have
\begin{align}
2^{kp^l}=\sum_{j=0}^{kp^l}\binom{kp^l}{j}
\equiv\sum_{j=0}^{kp^{l-1}}\binom{kp^{l-1}}{j}=2^{kp^{l-1}}
\pmod{p^l},\\
\binom{2kp^l}{kp^l}=\sum_{j=0}^{kp^l}\binom{kp^l}{j}^2
\equiv\sum_{j=0}^{kp^{l-1}}\binom{kp^{l-1}}{j}^2
=\binom{2kp^{l-1}}{kp^{l-1}}
\pmod{p^l},
\end{align}
found by using \eqref{eq:basic_binom_congruence} again.
Thus we obtain
\begin{equation}\label{eq:needed-2}
\binom{-\frac12}{kp^l}=(-1)^{kp^l}2^{-kp^l}\binom{2kp^l}{kp^l}
\equiv
(-1)^{kp^{l-1}}2^{-kp^{l-1}}\binom{2kp^{l-1}}{kp^{l-1}}
=\binom{-\frac12}{kp^{l-1}}
\pmod{p^l}.
\end{equation}
Combining \eqref{eq:needed-1}, \eqref{eq:needed-2} and
\eqref{eq:basic_binom_congruence},
we have
\begin{equation}
\binom{mp^r}{pj}\binom{-\frac12}{pj}\equiv
\binom{mp^{r-1}}{j}\binom{-\frac12}{j}\pmod{p^r}.
\end{equation}
Hence the desired $\J_2$-congruence relation,
\eqref{eq:higher_congruence_2}, follows.

Next we prove the $\J_3$-congruences.
Noting the relations
\begin{align*}
\J_3(p^r)p^{3r}
&=-2\sum_{k=0}^{p^r}(-1)^k\binom{-\frac12}{k}^{\!2}\binom{p^r}{k}
\sum_{j=0}^{k-1}\frac1{(2j+1)^3}\binom{-\frac12}{j}^{\!\!-2}p^{3r}\\
&\equiv-2\sum_{k=0}^{p^{r-1}}(-1)^k\binom{-\frac12}{k}^{\!2}\binom{p^{r-1}}{k}
\sum_{j=0}^{kp-1}\frac1{(2j+1)^3}\binom{-\frac12}{j}^{\!\!-2}p^{3r}\\
&\equiv-2\sum_{k=0}^{p^{r-1}}(-1)^k\binom{-\frac12}{k}^{\!2}\binom{p^{r-1}}{k}
\sum_{l=0}^{k-1}\frac1{(2l+1)^3}\binom{-\frac12}{\frac{p(2l+1)-1}2}^{\!\!-2}p^{3(r-1)},
\end{align*}
to prove \eqref{eq:higher_congruence_3}
it is enough to show the congruence 
\begin{equation}\label{eq:p-order_estimation}
\frac{p^{3(r-1)}}{(2l+1)^3}\left\{
\binom{-\frac12}{\frac{p(2l+1)-1}2}^{\!\!-2}
-\binom{-\frac12}{l}^{\!\!-2}
\right\}\equiv0
\pmod{p^r}
\end{equation}
for $0\le l\le p^{r-1}$.
It is immediately seen that
the validity of the congruence \eqref{eq:p-order_estimation}
is equivalent to
the validity of the inequality
\begin{equation}\label{eq:reducedlemma}
4\ord_p\binom{2j}{j}-\ord_p\left\{
\binom{p(2j+1)-1}{\frac{p(2j+1)-1}2}^2-\binom{2j}{j}^2
\right\}+\ord_p(2j+1)\le 2r-1
\end{equation}
for $0\le j\le p^r$.
Here we denote by $\ord_p(m)$ the exponent of $p$ in $m\in\N$
for a given prime number $p$.
This inequality is obtained as follows.
We first note that
\begin{equation}\label{eq:(6.15)}
\binom{mp^r-1}{k}\equiv(-1)^{k-[k/p]}\binom{mp^{r-1}}{[k/p]}
\pmod{p^r},
\end{equation}
where $[x]$ denotes the largest integer not exceeding $x$.
Then, utilizing \eqref{eq:(6.15)},
we have
\begin{align*}
\ord_p\binom{p(2j+1)-1}{\frac{p(2j+1)-1}2}&=\ord_p\binom{2j}{j},\\
\binom{p(2j+1)-1}{\frac{p(2j+1)-1}2}^2&\equiv\binom{2j}{j}^2
\pmod{p^{1+\ord_p(2j+1)}},
\end{align*}
which imply the equality
\begin{equation}
2\ord_p\binom{2j}{j}-\ord_p\left\{
\binom{p(2j+1)-1}{\frac{p(2j+1)-1}2}^2-\binom{2j}{j}^2
\right\}+\ord_p(2j+1)=-1
\end{equation}
for any $j\ge0$.
Finally, noting the obvious relation
$\ord_p\binom{2j}{j}\le r$ for $0\le j\le p^r$,
we obtain the desired result, \eqref{eq:higher_congruence_3}.
\end{proof}

\begin{rem}
While the Ap\'ery numbers
$A_n = \sum_{k=0}^n \binom{n}{k}^2\binom{n+k}{k}^2$
satisfy the supercongruence relation
\begin{equation}
A_{kp^n-1} \equiv A_{kp^{n-1}-1} \pmod{p^{3n}}
\end{equation}
for any prime $p$ and any $k\in\N$ (see, e.g., \cite{Be1987}),
it seems that no such relation exists for
the numbers $\J_2(n)$.
For instance, the congruence relation
\begin{equation}
\J_2(mp^r) \equiv \J_2(mp^{r-1}) \pmod{p^{r+1}}
\end{equation}
does not hold in general.
\end{rem}

\begin{rem}
Although we have not been able to obtain a proof, 
numerical computations indicate the relation
\begin{equation}
\sum_{k=0}^{p-1}\J_2(k)^2 \equiv \left(\frac{-1}p\right) \pmod{p^3}
\end{equation}
for any odd prime $p$.
This relation is quite similar to the
Rodriguez-Villegas-type congruence \cite{Mo2003},
\begin{equation}
\sum_{k=0}^{p-1}2^{-4k}\binom{2k}k^2 \equiv \left(\frac{-4}p\right) \pmod{p^2}.
\end{equation}
Here $\left(\frac{a}b\right)$ denotes the Legendre symbol. 
This similarity seems to suggest an algebra-geometric interpretation of
the numbers $\J_2(n)$ and $\J_3(n)$ and
special values $\zeta_Q(2)$ and $\zeta_Q(3)$.
\end{rem}

%====================================================
\section{Appendix}\label{sec:appendix}
%====================================================

Here we establish a general formula for the local holomorphic solution
of a particular inhomogeneous hypergeometric differential equation.

The (generalized) hypergeometric function is defined by
\begin{equation}
\HGF{p}q(\va;\vb;z)
\deq\sum_{n=0}^\infty\frac{(\va)_n}{(\vb)_n}
\frac{z^n}{n!}
\end{equation}
for $\va=(a_1,\dots,a_{p})$ and $\vb=(b_1,\dots,b_q)$.
Here, for simplicity, we introduce the following notations
\begin{equation*}
(a)_n\deq\dfrac{\Gamma(a+n)}{\Gamma(a)},\quad
(\va)_n\deq(a_1)_n\dotsb(a_{p})_n,\quad
(\vb)_n\deq(b_1)_n\dotsb(b_q)_n.
\end{equation*}
We assume that
the quantities $a_j$ and $b_j$ are not nonnegative integers.
The function $\HGF{p}q(\va;\vb;z)$ satisfies the differential equation
$\HGDop{p}q(\va;\vb;z)\cdot\HGF{p}q(\va;\vb;z)=0$,
where the operator $\HGDop{p}q(\va;\vb;z)$ is defined by
\begin{equation}
\begin{split}
\HGDop{p}q(\va;\vb;z)
&\deq-(\euler+a_1)\dotsb(\euler+a_{p})+\D_z(\euler+b_1-1)\dotsb(\euler+b_q-1),
\end{split}
\end{equation}
$\euler:=z\D_z$ being the Euler (degree) operator.
We first demonstrate the following.

\begin{lem}\label{prop:polynomial}
Let
\begin{equation}
P_d(\va;\vb;z)
\deq\sum_{n=0}^d \frac{(\va)_n}{(\vb)_n}\frac{z^n}{n!}.
\end{equation}
Then we have
\begin{equation}
\HGDop{p}q(\va;\vb;z)P_d(\va;\vb;z)
=-\frac{(\va)_{d+1}}{(\vb)_d}\frac{z^d}{d!}.
\end{equation}
\end{lem}

\begin{proof}
The assertion immediately follows from the equality
\begin{equation}
\HGDop{p}q(\va;\vb;z)\left(\frac{(\va)_n}{(\vb)_n}\frac{z^n}{n!}\right)
=-\frac{(\va)_{n+1}}{(\vb)_n}\frac{z^n}{n!}
+\frac{(\va)_n}{(\vb)_{n-1}}\frac{z^{n-1}}{(n-1)!}.
\end{equation}
\end{proof}

\begin{prop}\label{prop:inhomHGF}
Let $g(z)=\sum_{n=0}^\infty c_n z^n/n!$
be a holomorphic function around $z=0$.
Suppose that $\sum_{n=0}^\infty c_n\frac{(\vb)_n}{(\va)_{n+1}}$
converges absolutely.
Then the local holomorphic (power series at $z=0$) 
solution of the inhomogeneous
differential equation
\begin{equation}
\HGDop{p}q(\va;\vb;z)f(z)=g(z)
\end{equation}
with the initial condition $f(0)=C$ is given by
\begin{equation}
\begin{split}
f(z)&=\sum_{n=0}^\infty
\left(C+\sum_{k=0}^{n-1}c_k\frac{(\vb)_k}{(\va)_{k+1}}\right)
\frac{(\va)_n}{(\vb)_n}\frac{z^n}{n!}\\
&=C\cdot\HGF{p}q(\va;\vb;z)+\sum_{n=0}^\infty
\left(\sum_{k=0}^{n-1}c_k\frac{(\vb)_k}{(\va)_{k+1}}\right)
\frac{(\va)_n}{(\vb)_n}\frac{z^n}{n!}.
\end{split}
\end{equation}
\end{prop}

\begin{proof}
Write
\begin{equation}
F(z)=-\sum_{k=0}^\infty c_k\frac{(\vb)_k}{(\va)_{k+1}}P_k(\va;\vb;z).
\end{equation}
This series converges absolutely near the origin,
owing to the relation
\begin{equation*}
\sum_{k=0}^\infty \left|c_k\frac{(\vb)_k}{(\va)_{k+1}}P_k(\va;\vb;z)\right|
\le
\sum_{k=0}^\infty \left|c_k\frac{(\vb)_k}{(\va)_{k+1}}\right|
\sum_{n=0}^\infty \left|\frac{(\va)_n}{(\vb)_n}\frac{z^n}{n!}\right|.
\end{equation*}
It immediately follows that $\HGDop{p}q(\va;\vb;z)F(z)=g(z)$,
and hence $F(z)$ represents a local solution,
with $F(0)=-\sum_{k=0}^\infty c_k\frac{(\vb)_k}{(\va)_{k+1}}$.
Therefore, we find that
\begin{align*}
f(z)&\deq F(z)-F(0)\cdot\HGF{p}q(\va;\vb;z)\\
&=-\sum_{k=0}^\infty c_k\frac{(\vb)_k}{(\va)_{k+1}}\sum_{n=0}^k
\frac{(\va)_n}{(\vb)_n}\frac{z^n}{n!}
+\sum_{k=0}^\infty c_k\frac{(\vb)_k}{(\va)_{k+1}}\sum_{n=0}^\infty
\frac{(\va)_n}{(\vb)_n}\frac{z^n}{n!}\\
&=\sum_{n=0}^\infty\left(\sum_{k=0}^{n-1}c_k\frac{(\vb)_k}{(\va)_{k+1}}\right)
\frac{(\va)_n}{(\vb)_n}\frac{z^n}{n!}
\end{align*}
is the solution of $\HGDop{p}q(\va;\vb;z)f(z)=g(z)$
with the initial condition $f(0)=0$.
The assertion is now clear.
\end{proof}

As an application of this proposition,
we now give an alternative proof of Theorem \ref{thm:J3-formula}.
First, recall the definition of $\g_k(x)$ given in \S\ref{sec:preliminaries},
\begin{align}
\g_k(x)\deq \sum_{n=0}^\infty \binom{-\frac12}{n}\J_k(n)x^n.
\end{align}
Next we prove the following.

\begin{prop}
The functions $\g_2(x)$ and $\g_3(x)$ satisfy the differential equations
\begin{align}
\Wakayama \g_2(x) = 0,\qquad \Wakayama \g_3(x) = -\frac2{1+x},
\end{align}
where the operator $\Wakayama$ is given by
\begin{equation*}
\Wakayama = 8x^2(1+x)^2\D_x^3+24x(1+x)(1+2x)\D_x^2+2(4+27x+27x^2)\D_x+3(1+2x).
\end{equation*}
\end{prop}

\begin{proof}
Recall the relation
\begin{equation*}
\g_2(x)=\frac1{\sqrt{1+x}}\HGF32\left(\frac12,\frac12,\frac12;1,1;\frac{x}{1+x}\right).
\end{equation*}
It is easy to show that
the differential equation $\Wakayama \g_2(x)=0$ is equivalent to
\begin{equation*}
\HGDop32\left(\frac12,\frac12,\frac12;1,1;z\right)(\g_2(x)\sqrt{1+x})=0
\end{equation*}
for $z=\frac{x}{1+x}$.
In fact, we have
\begin{equation}
\frac18(1+x)^{\frac32}\Wakayama=
\HGDop32\left(\frac12,\frac12,\frac12;1,1;z\right)\sqrt{1+x}.
\end{equation}
Then, employing the recurrence relation \eqref{rec_for_J3} of $\J_3(n)$,
we have
\begin{align*}
\Wakayama \g_3(x) &=
\sum_{n=0}^\infty
\left\{
2n(2n+1)(2n-1)c_{n-1}+(2n+1)(8n^2+8n+3)c_n+8(n+1)^3c_{n+1}
\right\}x^n\\
&=-2\sum_{n=0}^\infty(-1)^nx^n=\frac{-2}{1+x},
\end{align*}
where $c_n=(-1)^n 2^{-2n}\binom{2n}{n}\J_3(n)$.
\end{proof}

Next, defining
$\varphi(z)\deq \g_3(x)\sqrt{1+x}$ for $z=\frac{x}{1+x}$,
we have
\begin{equation}
\begin{split}
\HGDop32\left(\frac12,\frac12,\frac12;1,1;z\right)\varphi(z)
&=\frac1{8(1-z)^{\frac32}}\Wakayama g_3(x)
=-\frac{1}{4\sqrt{1-z}}
=-\frac14\sum_{n=0}^\infty \left(\frac12\right)_{\!\!n} \frac{z^n}{n!}.
\end{split}
\end{equation}
Then, because $\varphi(0)=0$,
we find
\begin{equation}
\begin{split}
\varphi(z)&=
-\frac14\sum_{n=0}^\infty\left(\sum_{k=0}^{n-1}\left(\frac12\right)_{\!\!k}
\frac{(1,1)_k}{(\frac12,\frac12,\frac12)_{k+1}}\right)
\frac{(\frac12,\frac12,\frac12)_{n}}{(1,1)_n}\frac{z^n}{n!}\\
&=-2\sum_{n=0}^\infty\left(\sum_{k=0}^{n-1}\frac{1}{(2k+1)^3}
\binom{-\frac12}{k}^{\!\!\!-2}\right)
(-1)^n\binom{-\frac12}{n}^{\!\!3}z^n,
\end{split}
\end{equation}
using Proposition \ref{prop:inhomHGF}.
It follows that
\begin{equation}
\begin{split}
\g_3(x)&=
\frac{-2}{\sqrt{1+x}}\sum_{n=0}^\infty\left(\sum_{k=0}^{n-1}\frac{1}{(2k+1)^3}
\binom{-\frac12}{k}^{\!\!\!-2}\right)(-1)^n\binom{-\frac12}{n}^{\!\!3}
\left(\frac{x}{1+x}\right)^{\!n},
\end{split}
\end{equation}
from which we obtain Theorem \ref{thm:J3-formula}.

\smallskip

\begin{flushleft}
Kazufumi KIMOTO\\
Department of Mathematical Science,
University of the Ryukyus.\\
Senbaru, Nishihara, Okinawa 903-0231, JAPAN.\\
\texttt{kimoto@math.u-ryukyu.ac.jp}\\

\bigskip

Masato WAKAYAMA\\
Faculty of Mathematics,
Kyushu University.\\
Hakozaki, Fukuoka 812-8518, JAPAN.\\
\texttt{wakayama@math.kyushu-u.ac.jp}
\end{flushleft}

\end{document}